\documentclass[a4paper,12pt]{amsart}
\usepackage{amssymb,amsthm,bbm,epic,eepic,graphics,amsbsy,url}
\usepackage{a4wide}
\usepackage{amssymb,amsthm}
\usepackage{xypic}
\usepackage[applemac]{inputenc}
\usepackage{longtable}
\usepackage{array}
\newcommand{\la}{\lambda}

\DeclareMathOperator{\rk}{\mathrm{rk}}
\DeclareMathOperator{\Sp}{\mathrm{Sp}}
\DeclareMathOperator{\tr}{\mathrm{tr}}
\DeclareMathOperator{\Irr}{\mathrm{Irr}}

\def\Mat{\mathrm{Mat}}

\def\End{\mathrm{End}}
\def\Ker{\mathrm{Ker}\,}
\def\Res{\mathrm{Res}}

\def\Aut{\mathrm{Aut}}

\def\C{\ensuremath{\mathbbm{C}}}
\def\Q{\mathbbm{Q}}
\def\Z{\mathbbm{Z}}

\def\R{\mathbbm{R}}

\def\gl{\mathfrak{gl}}
\def\sl{\mathfrak{sl}}
\def\so{\mathfrak{so}}
\def\sp{\mathfrak{sp}}

\def\k{\mathbf{k}}
\def\kt{\mathbf{k}^{\times}}

\def\eps{\epsilon}

\def\onto{\twoheadrightarrow}
\def\into{\hookrightarrow}
\def\ii{\mathrm{i}}
\def\Gal{\mathrm{Gal}}
\def\Irr{\mathrm{Irr}}
\def\GL{\mathrm{GL}}
\def\SL{\mathrm{SL}}
\def\SO{\mathrm{SO}}
\def\SP{\mathrm{SP}}

\newtheorem{theo}{Theorem}[section]

\newtheorem{prop}[theo]{Proposition}

\newtheorem{lemma}[theo]{Lemma}

\newtheorem{cor}[theo]{Corollary}

\def\mod{\ \mathrm{mod}\ }

\makeatletter
\DeclareRobustCommand\widecheck[1]{{\mathpalette\@widecheck{#1}}}
\def\@widecheck#1#2{%
   \setbox\z@\hbox{\m@th$#1#2$}%
   \setbox\tw@\hbox{\m@th$#1%
      \widehat{%
         \vrule\@width\z@\@height\ht\z@
         \vrule\@height\z@\@width\wd\z@}$}%
   \dp\tw@-\ht\z@
   \@tempdima\ht\z@ \advance\@tempdima2\ht\tw@ \divide\@tempdima\thr@@
   \setbox\tw@\hbox{%
      \raise\@tempdima\hbox{\scalebox{1}[-1]{\lower\@tempdima\box\tw@}}}%
   {\ooalign{\box\tw@ \cr \box\z@}}}
\makeatother
\def\cnp#1#2{ \left( \begin{array}{c} {#1} \\ {#2} \end{array} \right)}
\title{{\bf Braids inside the Birman-Wenzl-Murakami algebra}}
\author{Ivan Marin}
\date{October 18, 2009}
\begin{document}

\maketitle

\bigskip
\begin{center}
Institut de Math\'ematiques de Jussieu \\
Universit\'e Paris 7 \\
175 rue du Chevaleret \\
F-75013 Paris
\end{center}
\bigskip

\bigskip

\bigskip

\noindent {\bf Abstract.} We determine the Zariski closure of
the representations of the braid groups that factorize through
the Birman-Wenzl-Murakami algebra, for generic values
of the parameters $\alpha,s$. For $\alpha,s$ of modulus 1 and
close to 1, we prove that these representations are unitarizable,
thus deducing the topological closure of the image when
in addition $\alpha,s$ are algebraically independent. 
\medskip


\def\RR{\mathcal{R}}
\def\SS{\mathcal{S}}
\def\HH{\mathcal{H}}
\def\BB{\mathcal{B}}
\def\XX{\mathrm{X}}

\section{Introduction}

Let $B_n$ denote the braid group on $n$ strands, defined by the presentation
with generators $\sigma_1,\dots, \sigma_{n-1}$ and relations
$\sigma_i \sigma_{i+1} \sigma_i = \sigma_{i+1} \sigma_i \sigma_{i+1}$,
$\sigma_i \sigma_j = \sigma_j \sigma_i$ for $| i - j | \geq 2$ -- which
imply that all generators are conjugated.
We
consider here the linear representations of $B_n$ afforded
by the so-called Birman-Wenzl-Murakami algebras. For a field
$K$ of characteristic 0 and $\alpha,s \in K$ with $\alpha,s,s-s^{-1}$ nonzero, the
algebra $BMW_n(s,\alpha)$ can be defined as the quotient of the group
algebra $K B_n$ by the 3 relations
$$
(\sigma_1 - s)(\sigma_1+s^{-1})(\sigma_1 + \alpha^{-1}) = 0
$$
and
$$
\left( 1- \frac{\sigma_2 - \sigma_2^{-1}}{s-s^{-1}} \right)
\sigma_1^{\pm 1} 
\left( 1- \frac{\sigma_2 - \sigma_2^{-1}}{s-s^{-1}} \right)
= \alpha^{\pm 1}
\left( 1- \frac{\sigma_2 - \sigma_2^{-1}}{s-s^{-1}} \right).
$$
For generic
values of $\alpha,s$, the algebra $BMW_n(s,\alpha)$ is semisimple
and its structure is known, thus providing many representations of the
braid groups.

The algebra $BMW_n(s,\alpha)$ is a deformation of Brauer's
centralizer algebra (see below), which admits for quotient
the so-called Iwahori-Hecke algebra of type $A_{n-1}$, namely
the quotient $H_n(s)$ of $K B_n$ by the relation
$$
(\sigma_1 - s)(\sigma_1 + s^{-1}) = 0.
$$
With this relation, rewritten as $\sigma_1 - \sigma_1^{-1} = s-s^{-1}$,
the last two relations of $BMW_n(s,\alpha)$ are void, making
$H_n(s)$ appear as a quotient of $BMW_n(s,\alpha)$. 

Besides the representations induced by this quotient, $BMW_n(s,\alpha)$
admits another special representation, already singled out
in \cite{BW}, which is known to induce a faithful
representation of $B_n$ by work of Krammer and Bigelow (see \cite{KRAM0,BIG,KRAM}). We call
it the Krammer representation $R_K : B \to \GL_{n(n-1)/2}(K)$.

Let $R : B_n \to \GL_N(K)$ be a linear representation afforded
by some linear representation of $BMW_n(s,\alpha)$. We are interested here
in the image $R(B_n) \subset \GL_N(K)$. Since $R$ is defined only
up to conjugacy, the first natural question is about the closure
of $R(B_n)$ for the Zariski topology.

Letting $B_n' = (B_n,B_n)$ denote the commutator subgroup, we
state some the results obtained in terms of $R(B_n')$, because
the statements are simpler. Since $B_n' = \Ker(B_n \onto \Z)$
is defined by $\sigma_1 \mapsto 1$, the corresponding
statements for $R(B_n)$ are easy to deduce from them.

For $\alpha,s$ generic (say, algebraically independant over $\Q \subset K$), this
problem was solved in \cite{LIETRANSP} for the Hecke algebra representations,
and in \cite{KRAMMINF} for the Krammer representation. The present paper is thus
a sequel of these two previous works. In \cite{LIETRANSP}
we proved that the Zariski closure of $B_n$ inside the whole Hecke
algebra had for Lie algebra the Lie subalgebra of the group algebra
of the symmetric group $\mathfrak{S}_n$ generated (for the bracket
$[a,b] = ab-ba$) by the transpositions, and decomposed this reductive
Lie algebra.
We called this Lie algebra the infinitesimal Hecke algebra (of type
$A_{n-1}$). Here we exhibit a reductive Lie subalgebra of the Brauer
centralizer algebra that plays a similar role, and that we decompose
accordingly. A consequence is the following, that generalize
\cite{KRAMMINF}, theorem 2. Recall from \cite{BW} that $BMW_n(s,\alpha)$
is split semisimple over $K$ (see \cite{BW} theorem 3.7, the
proof given there being valid over $\Q(s,\alpha)$, not only $\C(s,\alpha)$).

\begin{theo} \label{maintheo} Let $R : B_n \to \GL_N(K)$ an irreducible representation
afforded by $BMW_n(s,\alpha)$ for $\alpha,s$ algebraically independant over $\Q$
which does not factorize through $H_n(s)$. Then $R(B_n)$
is Zariski-dense in $\GL_N(K)$. If $R = R_0 \oplus R_1 \oplus \dots \oplus R_k
: B_n \to \GL_N(K)$
with $R_i : B_n \to \GL_{N_i}(K_i)$ as before for $i \geq 1$, $N = N_0+ N_1+\dots+N_k$,
$R_0$ factoring through $H_n(s)$,
and $R_i \not\simeq R_j$, then $R(B_n')$ is Zariski-dense in $G_0 \times \SL_{N_1}(K)
\times \dots \times \SL_{N_k}(K)$ where $G_0$ is the closure
of $R_0(B_n')$. The same holds if $\alpha = s^m$ for $m$ outside
a finite set of integer values.
\end{theo}

Another natural question is about the unitarisability of the
representations of $B_n$ obtained this way, when $K = \C$. In this case,
the determination of the Zariski closure done above is more or less
equivalent to the determination of the topological
closure for the usual topology of $\C$.

This question makes sense only for $|\alpha| = |s| = 1$, because
of the spectrum of the Artin generator. Even in this case,
it is known that the representations afforded by $H_n(s)$
are not unitarizable in general, but they are so if in addition
$s$ is close enough to 1. This was first proved by H. Wenzl, who
exhibited in \cite{WENZLHECKE} explicit unitary matrix models
for these representations.
For $\alpha,s$ close to 1 and some additional constraints,
this was also proved for the Krammer representation by R. Budney,
who followed the method of C. Squier (for the Burau representation) of constructing
an explicit sesquilinear form preserved in the usual matrix models
of this representation. 

We showed in previous works (see \cite{GT,KRAMMINF}) how to obtain new proofs of these two results by
making use of Drinfeld theory of associators. The idea is that
all these representations appear as monodromy of so-called
KZ-systems, and that these KZ-systems have for coefficients
real matrices which are compatible with certain natural bilinear forms. Then, from the choice of a Drinfeld associator
with rational (or real) coefficients, whose existence was proved
by Drinfeld, one can built a representation of $B_n$ over the ring $\R[[h]]$
of formal series with image in some formal unitary group
(with respect to the automorphism $f(h) \mapsto f(-h)$). By specializing
the matrices in $h$ a purely imaginary complex number we then get unitary
representations of $B_n$, and this provides a natural path for
explaining the unitarisability of this kind of representations.
A technical problem however is that we do not have (so far) any
insurance that the series involved have nonzero convergence radius.
For the Hecke algebra representations (see \cite{GT}), one
can explicitely compute the matrix models obtained this way,
which turn out to be the same than the ones obtained earlier by H. Wenzl, and check that
they are indeed convergent.
Another way, that we already used in \cite{KRAMMINF} for the Krammer
representation, is to approximate the representation over
$\R[[h]]$ by equivalent representations over the ring of convergent power
series. This is the device we use
here, giving in the above spirit a new proof of the following result,
which was originally proved by Wenzl (see \cite{WENZLBCD}) by using
the Jones construction.

\begin{theo} \label{theounit} (Wenzl) Let $S^1 = \{ z \in \C \ | \ |z| = 1 \}$. There exists
an open subset $U$ of $S^1 \times S^1$ whose closure contains
$(1,1)$ such that, for $(s,\alpha) \in U$, the representations
of $B_n$ induced by $BMW_n(s,\alpha)$ are unitarizable.
\end{theo}

Putting the two theorems together, they determine up to isomorphism the topological
closure of all representations of $B_n$ that factorize through $BMW_n(s,\alpha)$
for $(s,\alpha) \in U$
with $s,\alpha$ algebraically independant over $\Q$, since the
compact Zariski-dense subgroups of $\GL_N(\C)$ are its maximal compact
subgroups. This generalizes (part of) the work of Freedman, Larsen and Wang in \cite{FLW}
on the representations of $H_n(s)$.




Finally note that the open set $U$ of theorem \ref{theounit} contains
(a dense set of) couples $(s,\alpha)$ which are algebraically independent.
In particular, every irreducible representation $R$ as in theorem \ref{maintheo}
(not factorizing through $H_n(s)$) enables, when faithful, to embed
$B_n$ in the corresponding unitary group as a dense subgroup.

\section{Brauer diagrams and Brauer algebra}

We refer to \cite{WENZLBRAUER} or \cite{GOODMAN} for the definition
and basic properties of the algebra $Br_n$ of Brauer diagrams (in short : Brauer algebra)
as a finite-dimensional algebra over $\Q[m]$, and its specialization
$Br_n(m)$ over some field $\k$ of characteristic $0$, where
$m \in \k$. They are algebras spanned by so-called Brauer diagrams,
where $ab$ means composing the diagram $b$ below the diagram $a$
with additional relations $p_{ij}^2 = m p_{ij}$. These algebras
contain $\mathfrak{S}_n$, and in particular the transpositions $s_{ij}$.
In this section we assume $\k \subset \R$, and in particular $m \in \R$.

\begin{figure}
\begin{center}
\resizebox{!}{3cm}{
\includegraphics{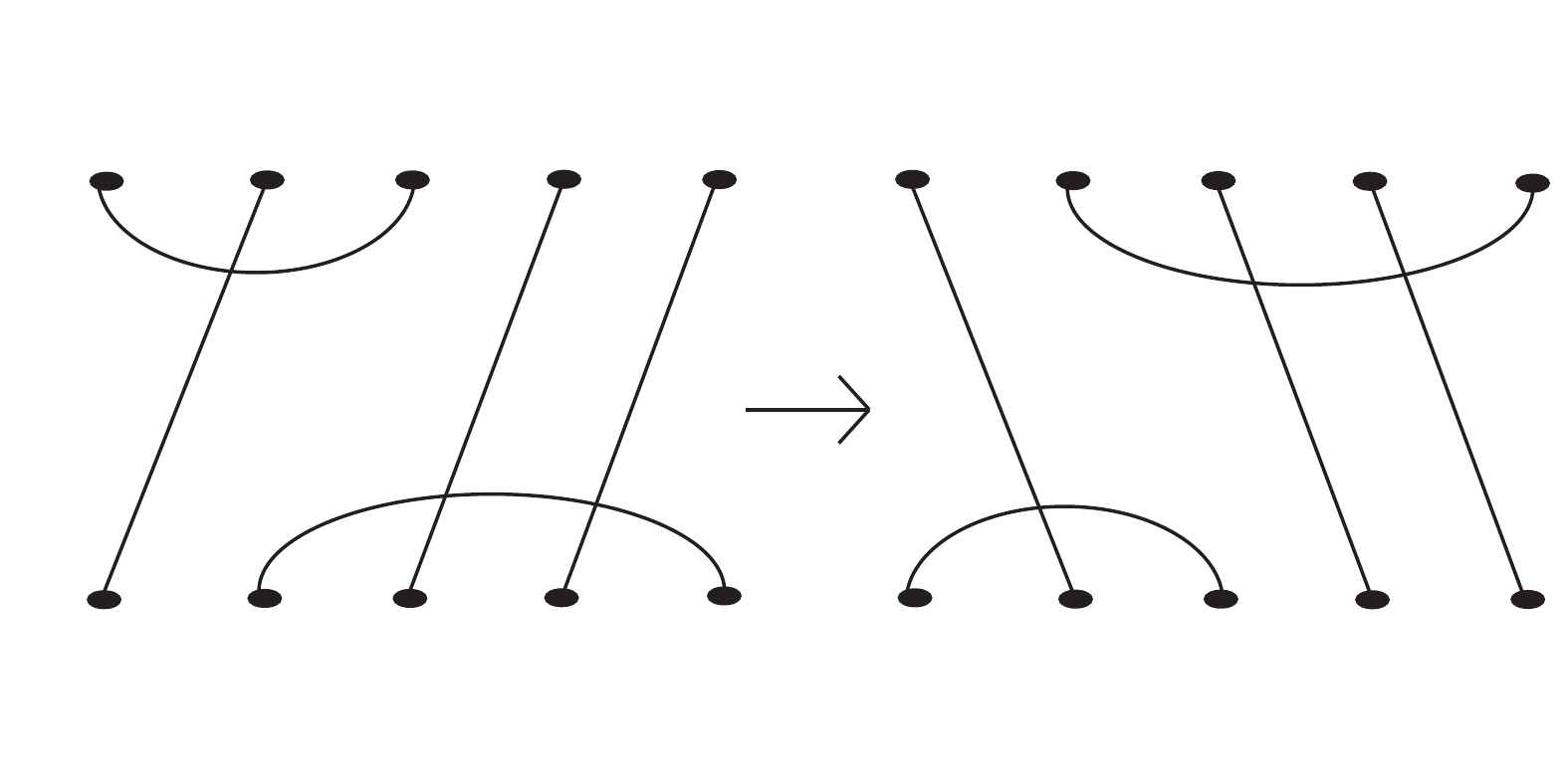}}
\end{center}
\caption{The `upside-down' automorphism $\tau$.}
\label{updown}
\end{figure}

We first define an involutive linear
automorphism $\tau$, defined at the level of the diagrams
by reflecting the diagram `upside-down' (see figure \ref{updown}).
For $w \in \mathfrak{S}_n$
we have $\tau(w) = w^{-1}$ ; moreover $\tau(p_{ij}) = p_{ij}$.
This anti-automorphism thus leaves the generators $s_{ij}$ and $p_{ij}$
invariant, and can alternatively be defined from this property. 
For $m > n$ or $m \not\in \Q$ there exists a non-degenerate trace
$\tr_M$ such that $\tr_M(b) = \tr_M(\tau(b))$ for every diagram $b$,
for instance the Markov trace defined in \cite{WENZLBRAUER}.
We assume from now on that $\tr_M$ is this Markov trace.
We then let $<D_1,D_2> = \tr_M(D_1 \tau(D_2))$ on the Brauer diagrams
and extend it by linearity. 
Consequences of the assumptions on $\tr_M$ are the following :
\begin{itemize}
\item For all $a,b$ we have
$$<a,b> = \tr_M(a \tau(b)) = \tr_M(\tau(a \tau(b))= \tr_M(b \tau(a)) = <b,a>.
$$
\item For all $a,b$, and $w \in \mathfrak{S}_n$,
$$
<wa,wb> = \tr_M(wa \tau(wb)) = \tr_M(wa \tau(b) \tau(w)) = \tr_M(\tau(w)
w a \tau(b)) = \tr_M(a \tau(b)) = <a,b>
$$
\item For all $a,b$, and $p_{ij}$,
$$
<p_{ij}a,b> = \tr_M(p_{ij}a \tau(b)) = \tr_M(a \tau(b) \tau(p_{ij})) = 
\tr_M(
 a \tau(p_{ij} b)) =  <a,p_{ij}b>
$$
\end{itemize}

In particular, the endomorphism $s_{ij} - p_{ij}$ is selfadjoint
with respect to $<\ , \ >$ and the elements of $\mathfrak{S}_n$
act orthogonally. For $\tr_M$ the Markov trace of \cite{WENZLBRAUER},
the bilinear form $(a,b) \mapsto \tr_M(ab)$ on $Br_n(m)$ is nondegenerate
for $m > n$ or $m \not\in \Q$. The Brauer algebra is then
semisimple and decomposes as a sum of matrix algebras. The set
$\Irr_n$ of irreducible representations of $Br_n$ is
in 1-1 correspondance with the partitions of $r$ with $n-r$ a
nonnegative even integer. When no confusion can arise, we denote them
by the corresponding partition $\la = (\la_1 \geq \la_2 \geq \dots)$,
or by $\la_n$ if the number of strands is not implicit. We
let $|\la| = \la_1 + \la_2 + \dots = r \leq n$.

For $\la \in \Irr_n$, let $\tr_{\la} : x \mapsto \tr(\la(x))$ denote the matrix
trace on the corresponding factor of $Br_n(m)$. By \cite{WENZLBRAUER},
we have
$$
\tr_M(b) = \sum_{\la \in \Irr_n} \frac{P_{\la}(m)}{m^n} \tr_{\la}(b).
$$
for some rational polynomials $P_{\la}(m)$ of $m$.
For $m > n$ it is possible to choose
the representations $\la \in \Irr_n$ defined over $\R$ and such that $\la(\tau(b)) = \ ^t \la(b)$.
Indeed, this means that the generators $s_{k,k+1}$ and $p_{12}$
have real entries and are symmetric in some basis, and it possible to find
such a basis by \cite{NAZAROV} (3.11) (see the proof of theorem 3.12).
Actually, the argument in \cite{NAZAROV} proves this under the
additional condition that $m$ is an integer. This additional condition
can be readily dropped, as the symmetric formulas obtained there
make sense for a real number $m > n$ hence define representations
of the corresponding $Br_n(m)$, the defining relations of $Br_n(m)$
being polynomial in $m$ and the matrix entries being (square
roots of) rational fractions of $m$.

We then have
$$
<a,b> = \sum_{\la \in \Irr_n} \frac{P_{\la}(m)}{m^n} \tr(\la(a \tau(b)))
= \sum_{\la \in \Irr_n} \frac{P_{\la}(m)}{m^n} \tr(\la(a ) \la(\tau(b)))
= \sum_{\la \in \Irr_n} \frac{P_{\la}(m)}{m^n} \tr(\la(a ) ^t\la(b))
$$

It follows that $<\ ,\ >$ is positive definite, for $m \not\in \Q$ or $m > n$,
if and only if all the $P_{\la}(m)/m^n$ are positive. From
\cite{WENZLBRAUER} we have an explicit a combinatorial description
of $P_{\la}(m)$, and we know that
$P_{\la}(m)$ coincides with the dimension of a representation
for $m$ a sufficiently large integer. In particular,
$P_{\la}(m) > 0$ for $m\gg 0$.
We let $\mathcal{S}_n  = \{ m \ | \ \exists \la \vdash r \leq n \ \ 
P_{\la}(m) = 0 \} \subset \Z$.

\section{Convergent approximation}

In this section we prove a technical result concerning
finitely generated subfields of the field $\R(\{h \})$ of convergent
Laurent series (i.e. the field of fractions of the ring $\R \{\{ h \} \}$ of
formal power series with nonzero convergence radius),
and apply it to the twisting of representations by adequate field
automorphisms. We denote
$\R((h))$ the field of formal Laurent series, $\R(h)$ the field of rational
fractions.

\begin{prop} Let $K = \R((h))$, $K^* = \R(\{ h \})$, $\eps \in \Aut(K)$
defined by $f(h) \mapsto f(-h)$, and $L \subset K$ a
finitely generated extension of $\R(h)$. There exists
$L^* \subset K^*$ 
such that $L^* \simeq L$ and an isomorphism $\Omega : L \to L^*$ with
$L \cap K^* \subset L^*$ and $\Omega(L \cap \R[[h]]) \subset L^* \cap \R\{\{h \} \}$. If $\eps(L) = L$, then $L^*$ can
be chosen such that $\eps(L^*) = L^*$ and the isomorphism $\Omega$ 
can be chosen such that $\eps \circ \Omega = \Omega \circ \eps$. Morever,
for any finite set $\{ u_1,\dots,u_t \} \subset L$ and $M > 0$,
$\Omega$ can be chosen such that $\Omega(u_i) \equiv u_i$ modulo $h^M$.
\end{prop}

\begin{proof}
Let $L_0 = L \cap K^*$ and $L_1$ a purely transcendantal extension of
$L_0$ in $L$. Since $L$ is finitely generated, $L/L_0$ has finite
transcendance degree $r$, and $L/L_1$ is finite, that it $L = L_1[\alpha]$
for some $\alpha \in L$ algebraic over $L_1$. Since $L_0 \supset
\R(h)$ we can assume $\alpha \in h \R[[h]]$ and also choose a trancendance
basis $f_1,\dots,f_r$ of $L_1/L_0$ with $f_i \in \R[[h]]$.
Since $L_0$ is finitely generated and $\R(\{ h \})$ has infinite transcendance
degree over $\R(h)$ (consider e.g. the family $\exp(h^d), d \geq 0$),
there exists $g_1,\dots g_r \in \R\{ \{ h \} \}$ algebraically
independant over $L_0$. Since $L_0 \supset \R(h)$, for any $P_1,\dots,
P_r \in \R[h]$, the family $g_1+P_1+\dots+g_r + P_r$ is also
algebraically independant over $L_0$, hence for any given $N$
we can choose $g_i \equiv f_i$ modulo $h^N$. Let then
$L_1^* = L_0(g_1,\dots,g_r) \simeq L_1$ and $P \in L_1[X]$ a minimal
polynomial for $\alpha \in L$. By not requiring $P$ to be monic,
since $L_0 \supset \R(h)$ we can assume that the coefficients of $P$
belong to $L_0[f_1,\dots,f_r] \cap \R[[h]]$. Through $L_1^* \simeq L_1$,
we define $P_g \in L_0[g_1,\dots,g_r]$. We have $P(\alpha) = 0$,
$P'(\alpha) \equiv \beta h^s$ modulo $h^{s+1}$ for some $s \geq 0$
and $\beta \in \R^{\times}$. By choosing $N$ large enough, we get $P'_g(\alpha)
\equiv \beta h^s$ modulo $h^{s+1}$ and $P_g(\alpha)\equiv 0$
modulo $h^{2s+1}$. In particuler $P_g(\alpha) \in P'_g(\alpha)^2 h \R[[h]]$.
Hensel's lemma then asserts (see e.g. \cite{EISENBUD} theorem 7.3)
the existence of $\gamma \in \R[[h]]$ such that $P_g(\gamma) = 0$
and $\gamma-\alpha \in P'_g(\alpha)h\R[[h]] \subset h \R[[h]]$. It follows
that $\gamma \in h \R[[h]]$. Then M. Artin's approximation theorem
(see \cite{MARTIN})
states that there exists a root $\tilde{\gamma} \in \R\{\{ h \} \}$ of $P_g$
that can be chosen arbitrarily close to $\gamma$, hence $\gamma \in \R \{ \{ h \} \}$
(since $P_g$ admits finitely many roots), and $L^* = L_1[\gamma] \simeq L$
is a subfield of $K^*$ that satisfies the wanted properties.

For any finite subset $u_1,\dots,u_t \in L = L_0(f_1,\dots,f_r)[\alpha]$,
by choosing $N$ large enough we can assume that $s \geq M$, hence
$\gamma \equiv \alpha$ modulo $h^M$,
and that $\Omega(u_i ) \equiv u_i$ modulo $h^M$.

Finally, assume that $\eps(L) = L$ and let $L^{\eps}
= \{ f \in L \ | \ \eps(f) = f \}$.
We have $\R(h^2) \subset L^{\eps} \subset K^{\eps} = \R((h^2))$,
$L = L^{\eps} \oplus h L^{\eps} \simeq L^{\eps}[X]/(X^2 -h)$. Let
$\Phi : K^{\eps} \to K$ be the isomorphism
defined by $f(h^2) \mapsto f(h)$ and $\Lambda = \Phi(L^{\eps})$.
Clearly $\Phi$ and $\Phi^{-1}$ map convergent series to convergent series.
We have $L \simeq \Lambda_+ = \Lambda[X]/(X^2-h)$ and
$\Gal(\Lambda_+/\Lambda) \simeq \Z/2$ is generated by the
action of $\eps$ on
$L \simeq \Lambda_+$. We have already show how to construct
an extension $\R(h) \subset \Lambda^* \subset K^*$ with $\Lambda^* \simeq
\Lambda$ over $\Lambda \cap K^* \supset \R(h)$. Therefore, there
exists an isomorphism between the extensions $\Lambda_+/\Lambda$ and
$\Lambda_+^*/\Lambda^*$ with $\Lambda_+^* = \Lambda[X]/(X^2-h)$.
Letting $L_-^* = \Phi^{-1}(\Lambda^*)$ and $L^* =
L_1^*(h) \subset K^*$ this defines $\Omega : L \to L^*$ with
$\Omega \circ \eps = \eps \circ \Omega$. Moreover, if $f \in
L \cap K^*$, then $f = f_1 + h f_2$ with $f_i \in L^{\eps} \cap
K^*$, hence $\Phi(f_i) \in \Lambda \cap K^*$, $\Phi(f_i) = f_i$ and
$\Omega(f) = f$. The verication that $\Omega$ can be chosen with $\Omega(u_i)$
close to $u_i$ in this case too, is straightforward and left to the reader.
\end{proof}

For $L$ a subfield of $\R((h))$ such that $\eps(L) = L$,
we let $U_N^{\eps}(L) = \{ X \in \GL_N(L) \ | \ X^{-1} = ^t \eps(X) \}$.

\begin{cor} Let $G$ be a finitely generated group, $R : G \to \GL_N(\R[[h]])$
be a linear representation, $g_1,\dots,g_r \in G$ such that
$R(g_1),\dots,R(g_r)$ have entries in $\R\{\{ h \}\}$, and $M > 0$.
Then there exists a linear representation $R^* : G \to \GL_N(\R\{ \{ h\} \}$
such that $R^*(g_i) = R(g_i)$ and, for all $g \in G$, $R^*(g) \equiv R(g)$ modulo $h^M$.
Moreover, if $R(G) \subset U_N^{\eps}(\R((h)))$, then we can assume
$R^*(G) \subset U_N^{\eps}(\R(\{ h \}))$.
\end{cor}
\begin{proof} Let $u_1,\dots,u_n$ be generators of $G$, and $L$
the subfield of $\R((h))$ generated over $\R(h)$ by the entries
of the $R(u_i)$ and $R(u_i^{-1})$. We have $R(G) \subset \GL_N(L)$
and $L$ is finitely generated over $\R(h)$.
Using an isomorphism $\Omega : L \to L^* \subset  \R(\{ h \})$
provided by the proposition, we let $R^* = \Omega \circ R$. If
$\Omega$ is chosen so that $R^*(u_i) \equiv R(u_i) \mod h^M$ and
$R^*(u_i^{-1}) \equiv R(u_i^{-1}) \mod h^M$ then $R^*(g) \equiv R(g) \mod h^M$
for every $g \in G$, which concludes the proof. In case $R(G) \subset
U_N^{\eps}(\R((h)))$, then $L$ clearly satisfies $\eps(L) = L$
and the condition $\Omega \circ \eps = \eps \circ \Omega$
implies $R^*(G) \subset U_N^{\eps}(\R(\{ h \}))$.
\end{proof}

\begin{prop} In the situation of the corollary, if $R$ is
absolutely irreducible, then $R^*$ can be chosen absolutely
irreducible.
\end{prop}
\begin{proof} Since $R$ is absolutely irreducible,
the image of the group algebra $\R((h)) G$ inside $Mat_N(\R((h)))$ is full, that
is there exists $g_1,\dots,g_{N^2} \in \R[[h]] G$ such that
$R(g_1),\dots,R(g_{N^2})$ is linearly independant over $\R((h))$.
The determinant of this family is thus congruent to $\beta h^s$
modulo $h^{s+1}$, for some $\beta \in \R^{\times}$ and $s \geq 0$.
Replacing the coefficients of $g_1,\dots,g_{N^2}$ by
their approximation modulo $h^{s+1}$ we can assume that
$g_1,\dots,g_{N^2}$ have coefficients in $\R[h]$.
Let $u_1,\dots,u_r$ denote the elements of $G$ that appear
in the $g_1,\dots,g_r$. Choosing $R^*$ such that $R^*(u_i) \equiv
R(u_i)$ modulo $h^{s+1}$ we get that $R^*(g_i) \equiv R(g_i) \mod h^{s+1}$,
whence the family $R^*(g_1),\dots,R^*(g_N)$ is also linearly
independant over $\R(\{ h \})$, which concludes the proof.
\end{proof}

\section{Unitarisability of the representations of the BMW algebras}

We first construct representations of the braid groups over
$\R[[h]]$ which are formally unitary, then approximate
these by convergent series. By an adequate specialization
this affords unitary representations which are shown to be equivalent
to representations of the BMW algebras.

Recall from e.g. \cite{DRINFELD} that the Lie algebra of (pure) infinitesimal braids
$\mathcal{T}_n$, or horizontal chord diagrams, is defined by generators
$t_{ij}$, $1 \leq i,j \leq n$ with relations $t_{ii} = 0, t_{ij} = t_{ji}$,
$[t_{ij}, t_{ik} + t_{kj}] = 0$, and $[t_{ij},t_{kl}] = 0$
for $\# \{ i,j,k,l \} = 4$. It is endowed with an action of $\mathfrak{S}_n$
by $w.t_{ij} = t_{w(i),w(j)}$, and a grading given by $\deg t_{ij} = 1$.
Drinfeld theory of associators (see \cite{DRINFELD}) defines, for $\k$ a field of
characteristic 0 and $\mu \in \kt$, a set $M_{\mu}(\k)$ of formal
series in two non-commuting variables, such that any
$\Phi \in M_{\mu}(\k)$ provides an (injective) morphism
$B_n \to \mathfrak{S}_n \ltimes \exp \widehat{\mathcal{T}}_n$
where $\widehat{\mathcal{T}}_n$ denotes the completion
of $\mathcal{T}$ with respect to its natural grading. Such a morphism
maps $\sigma_1$ to $(1 \ 2) \exp(\mu t_{12})$ (our $M_{\mu}(\k)$
is Drinfeld's $M_{2 \mu}(\k)$). Drinfeld then exhibits a special element
$\Phi_{KZ} \in M_{\ii \pi}(\C)$ and deduces abstracly from this that
$M_{\mu}(\k) \neq \emptyset$ for every $\mu$ and $\k$.
We refer to \cite{DRINFELD} for the main properties of such elements.

Let $\rho : \mathfrak{S}_n \to \GL_N(\R)$
be a representation of the symmetric group $\mathfrak{S}_n$
and $\varphi : \mathcal{T}_n \to \gl_N(\R)$ be a representation
of the Lie algebra of infinitesimal braids compatible with $\rho$,
i.e. $\rho(w) \varphi(t_{ij}) \rho(w)^{-1} = \varphi(t_{w(i),w(j)})$.
We can extend $\varphi$ into a representation $\widehat{\mathcal{T}}_n \to
\gl_N(\R[[h]])$ through $t_{ij} \mapsto h \varphi(t_{ij})$.
The choice of a real Drinfeld associator $\Phi \in M_1(\R)$ provides
a representation $R : B_n \to \GL_N(\R[[h]])$.

\begin{prop} If $\R^N$ is endowed with its canonical euclidean structure and
$\rho(\mathfrak{S}_n) \subset O_N(\R)$, $^t\varphi(t_{ij}) = \varphi(t_{ij})$,
then $R(B_n) \subset U_N^{\eps}(K)$.
\end{prop}
\begin{proof}
We only need to check that, for every Artin generator
$\sigma_i$, we have $R(\sigma_i) \subset U_N^{\eps}(K)$.
Recall from \cite{DRINFELD}
that $R(\sigma_i)$ is conjugate to $\exp h \varphi(t_{i,i+1}) \in
U_N^{\eps}(K)$ by some element of the form $\Phi(hx,hy)$
where $x,y$ are linear combination of the $\varphi(t_{ij})$, hence
are symmetric matrices. Now $\Phi$ is the exponential of a Lie series $\Psi$,
and $u \mapsto - ^t u$ is an automorphism of $\gl_N(\R)$. We thus get
$- ^t\Psi(hx,hy)= \Psi(- ^t h x, - ^t h y) = \eps(\Psi(hx,hy))$,
hence $ ^t \Phi(hx,hy)^{-1} = \eps( \Phi(hx,hy))$. It follows
that $\Phi(hx,hy) \in U_N^{\eps}(K)$ whence
$R(B_n) \subset U_N^{\eps}(K)$.
\end{proof}

We also recall from \cite{LIETRANSP} the following.

\begin{prop} \label{irred} If $\varphi$ is absolutely irreducible, then $R$ is
absolutely irreducible.
Moreover, the Lie algebra of the Zariski closure of $R(B_n)$
inside $\GL_N(K)$ contains $\varphi(\mathcal{T}) \otimes_{\R} K$.
\end{prop}
\begin{proof}
The proof of the first part follows from the fact that $R(\gamma_{ij})$
is congruent to $1 + 2h \varphi(t_{ij})$ modulo $h^2$, where $\gamma_{ij}$
denote the standard generator of the pure braid group, and that
$\varphi(\mathsf{U} \mathcal{T})$ is generated by the $\varphi(t_{ij})$.
By Burnside theorem there exists a basis of $\Mat_N(\R) = \R^{N^2}$ of
non-commutative polynomials $P_1,\dots,P_{N^2}$
in the $\varphi(t_{ij})$, that is $P_k = P_k( (\varphi(t_{ij})_{i,j})$.
Then the $P_k( ((\gamma_{ij}-1)/2h)_{ij})$ define elements of
$K B_n$ whose image under $R$ is congruent modulo $h$ to a basis
of $\Mat_N(\R)$. It follows that $R(KB_n)$ generates $\Mat_N(K)$,
that is that $R$ is absolutely irreducible. The last assertion
is an elementary consequence of Chevalley's formal exponentiation
theory, and is proved in \cite{LIETRANSP}, lemme 21.
\end{proof}

We then approximate $R$ by a convergent $R^* : B_n \to \GL_N(K^*)$
as in the previous section. A first remark is that $R^*$
can be chosen such that $R^*(\sigma_i)$ is conjugate to
$R(\sigma_i)$ in $\GL_N(K)$, at least when $R(\sigma_i)$
is diagonalisable with eigenvalues in $K^*$.
Indeed, if this is the case, there exists a polynomial $P \in K^*[X]$
with simple roots in $K^*$ such that $P(R(\sigma_i)) = 0$,
hence $P(R^*(\sigma_i) = \Omega(P) (\Omega(R(\sigma_i))) = \Omega(0) = 0$.
Moreover, the traces of the $R(\sigma_i)^k$, which belong to $K^*$,
equal the traces of the $R^*(\sigma_i)^k$, hence $R(\sigma_i)$
and $R^*(\sigma_i)$ have the same spectrum with multiplicities.
Notice that this assumption is satisfied in our case as soon as
$\varphi(t_{12})$ is semisimple, since $R(\sigma_1)=\rho(s_{12})
\exp(h \varphi(t_{12}))$.

Let $\alpha > 0$ such that the entries of the $R(\sigma_i)$
and $R(\sigma_i^{-1})$ all have convergence radius at least $\alpha$.
By specialization of $h$ to a purely imaginary number $\ii u$ of modulus
less than $\alpha$, we get unitary representations $R_{\ii u} : B \to U_N
\subset \GL_N(\C)$.

We now apply this to an irreducible representation of the Brauer
algebra $\la \in \Irr_n$, given in matrix form over $\R$ such that
$\la(\tau(b)) = ^t \la(b)$ for all $b$. Recall from \S 2 that this
is possible for $m > n$.
Then $\rho$ restricts to an orthogonal representation of $\mathfrak{S}_n$,
and $\varphi(t_{ij}) = \rho(s_{ij}) - \rho(p_{ij})$ defines a
compatible representation of $\mathcal{T}_n$ with $ ^t \varphi(t_{ij}) = \varphi(t_{ij})$.
We thus get representations
$R$ and $R^*$ of $B_n$.

\begin{prop} Assume $m \not\in \Q$ or $m > n$. For $m$ outside a finite
set of other values, the representations $R$ and $R^*$ factorize through
the Birman-Wenzl-Murakami algebra $BMW_n(e^h,e^{(1-m)h})$ and correspond to the same
partition $\la$.
\end{prop}
\begin{proof}
The assertion about $R$ is well-known, and can be proved e.g. along
the lines of \cite{QUOT}, proposition 4. Since the
defining relations of BMW have coefficients in $K^*$, then
$R^*$ also factorizes through $BMW$. Introduce
the elements $\delta_2= \sigma_1^2,\delta_3 = \sigma_2 \sigma_1^2 \sigma_2,\dots,
\delta_{n-1} = \sigma_{n-1}\dots \sigma_2 \sigma_1^2 \sigma_2 \dots \sigma_{n-1}$.
We have $R(\delta_k) = \exp h \varphi(Y_k)$, where
$Y_2 = t_{12}, \dots, Y_n = t_{1n} + \dots + t_{n-1,n}$.
We recall that, for $m > n$ or $m \not\in \Q$, and possibly
outside a finite set of other values of $m$, the values of
$\varphi(Y_2),\dots,\varphi(Y_n)$ determine $\la$ (see \cite{QUOT} \S 9).
Since the $R(\delta_k)$ have entries in $K^*$, we have
$R^*(\delta_k) = R(\delta_k)$ and the conclusion follows.
\end{proof}

\begin{figure}
\begin{center}
\resizebox{!}{3cm}{
\includegraphics{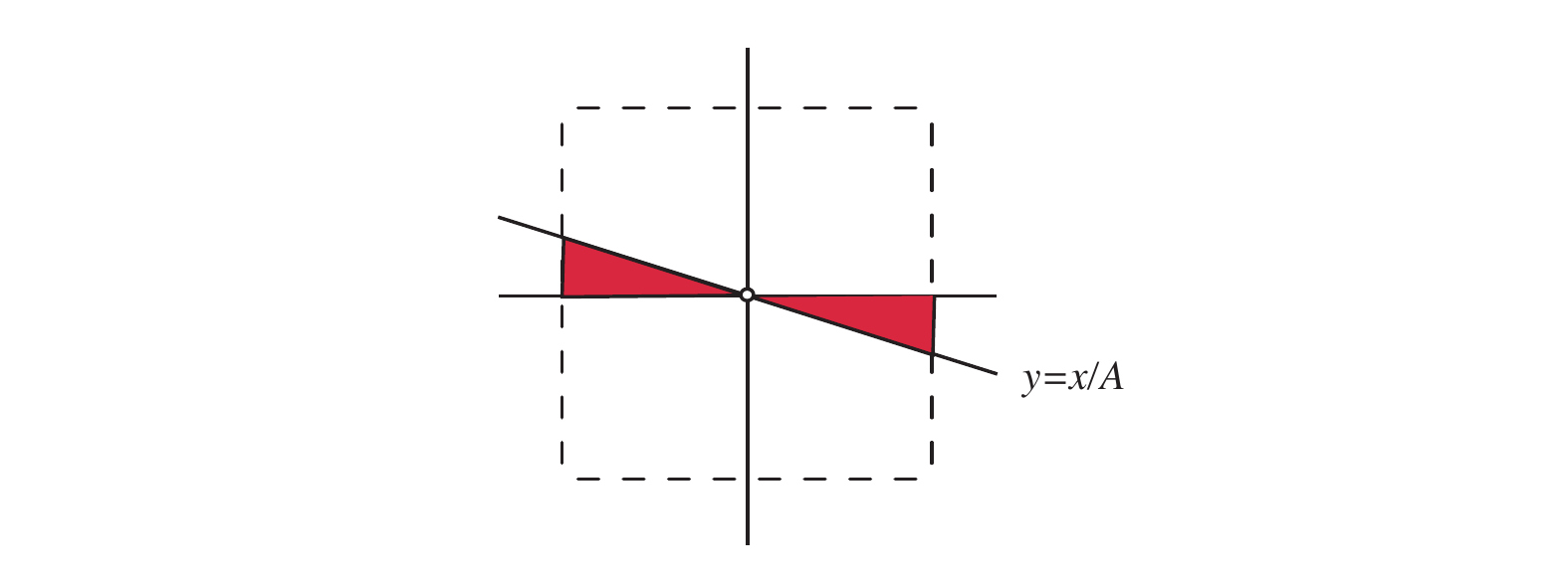}}
\end{center}
\caption{ $s = e^{\ii \pi y}, \alpha = e^{\ii \pi x}$}
\label{figS1S1}
\end{figure}

We denote $\log z$ the determination of the complex logarithm over
$\C \setminus \R_-$ such that $\log 1 = 0$. The following corollary
proves theorem \ref{theounit}, for $U$ drawn in figure \ref{figS1S1}.

\begin{cor} Assume $s, \alpha \in \C$ with $|\alpha| = |s| = 1$,
and $s, \alpha \not\in \{ 1 ,-1 \}$. There exists $A < -1$ and
$0 < \eta < 2$ such that, for $0 < |s-1| < \eta$, $0 < |\alpha-1|< \eta$
and $\frac{\log \alpha}{\log s }< A$, all the representations
of $B_n$ originating from $BMW_n(s,\alpha)$ are unitarizable.
\end{cor}

\begin{proof}

Let $m \in \R$ defined by $1-m = (\log \alpha)/(\log s)$. We can
choose $A$ and $\eta_0$ such that $m>n$ and does not belong to the
exceptions of the proposition for any $\la \in \Irr_n$, and such that,
for $0 < |s-1| < \eta_0$ and $0 < |\alpha-1| < \eta_0$, the algebra
$BMW_n(s,\alpha)$ is split semisimple and has for irreducible
representations the ones deduced from the generic algebras $BMW_n$
over $\R(\tilde{s},\tilde{\alpha})$, where $\tilde{s}$ and $\tilde{\alpha}$
are formal parameters. We mean by that that the
given irreducible representation $R_0$ of $B_n$ over $\C$ factorizing
through $BMW_n(s,\alpha)$ is deduced from the corresponding
representation $R_g$ of $B_n$ over $\R(\tilde{s},\tilde{\alpha})$ by
specialization of matrices, where the smallness of $\eta_0$ ensures that
$s,\alpha$, lie outside the poles of the entries.

The representation $R^*$ of $B_n$ over $K^* = \R(\{ h \})$ afforded
by the proposition satisfies $R^* \simeq R_g^*$ over $K^*$
where $R_g^*$ denotes the specialization of $R_g$ to
$\tilde{s} = e^h$ and $\tilde{\alpha} = e^{(1-m)h}$ (here again
we can change $A$ so that $e^h$, $e^{(1-m)h}$ lie outside
the poles, since $e^h$ and
$e^{(1-m)h}$ are algebraically independant over $\R$ for $m \not\in \Q$).
It follows that there exists $P \in \GL_N(K^*)$ such that
$R^*(b) = P R_g^*(b)P^{-1}$ for all $b \in B_n$. Let now $\beta > 0$
such that, for $0 < |h| < \beta$, the entries
of the matrices $R^*(\sigma_i)$,
$R_g(\sigma_i)$, $P$ and $P^{-1}$ are convergent. We choose $\eta$
such that $0 < \eta < \eta_0$ and $|e^{\ii u} -1| < \eta \Rightarrow
|u| < \beta$ for $u \in ]-\pi,\pi[$. Then $R_0$ is isomorphic to $R_{\ii u}^*$
which is a unitary representation of $B_n$. Since $\Irr_n$ is
finite we can choose $A,\eta$ uniformly in $\la$ and this concludes the
proof.
\end{proof}


\section{Infinitesimal Brauer algebras}

Let $\k$ be a field of characteristic 0, $m \in \k$ and $Br_n(m)$
be defined over $\k$.
In this section we study the Lie subalgebra of
$Br_n(m)$, with bracket given by $[a,b] = ab-ba$, which
is generated by the elements $t_{ij} = s_{ij} - p_{ij}$, that is the
images of the generators of $\mathcal{T}_n$ also denoted $t_{ij}$ in \S 4. We call
it the infinitesimal Brauer algebra and denote it $\BB_n(m)$.
The purpose of this section is to show that it is a reductive Lie algebra
for generic values of $m$ and to determine its structure in this case.
We recall from \S 2 that $\mathcal{S}_n$ denotes the set of values of $m$ for
which $Br_n(m)$ is not `generic'. We have $\mathcal{S}_n \subset \Z\cap
]-\infty,n]$.

\subsection{Reductiveness and Center}

We let $t = t_{12} = s_{12} - p_{12}$ and $s = s_{12}$.
A straightforward computation in $Br_2(m)$ shows that
$t^2 -1 = (m-2) p_{12}$ and $(t^2 - 1)(t+(m-1))=0$.
It follows that, for $m \neq 2$, $s_{ij}$ and $p_{ij}$ are polynomials
of $t_{ij}$. Since these elements generate $Br_n(m)$ as an algebra,
it follows that every irreducible representation $\rho$ of $Br_n(m)$
induces an irreducible representation $\rho_{\BB}$ of $\BB_n(m)$. In
particular, when $m \not\in \mathcal{S}_n$, $\BB_n(m)$ is a
reductive Lie algebra, as it admits a faithful semisimple
representation (take the direct sum of all irreducible representations
of $Br_n(m)$).

For $m \not\in \mathcal{S}_n$, $Br_n(m)$ is reductive as a Lie algebra.
We denote $p : Br_n(m) \onto Z(Br_n(m))$ the natural projection,
with kernel $[Br_n(m),Br_n(m)]$. Let $\mathbbm{T}$ denote the subspace
of $Br_n(m)$ spanned by the $t_{ij}$, and $T = \sum t_{ij} \in \mathbbm{T}$.
We have $T \in Z(\BB_n(m)) \subset Z(Br_n(m))$ for $m \neq 2$. Let
$\mathbbm{T}'$ denote the subspace of $\mathbbm{T}$ spanned by the
$t_{ij} - t_{kl}$. Since the $t_{ij}$ are linearly independent
we have $
\mathbbm{T} = \k T \oplus \mathbbm{T}'$. Since $\mathfrak{S}_n$
is 2-transitive, there exists $\sigma \in \mathfrak{S}_n$ with
$t_{kl} = t_{\sigma(i) \sigma(j)} = \sigma t_{ij} \sigma^{-1}$.
Now an explicit formula for $p$ can be given in terms
of the primitive idempotents $j_{\la}$ of $Br_n(m)$ for $\la \in \Irr_n$,
as $p(x) = \sum d(\la) tr_{\la}(x) j_{\la}$ where $d(\la)$
is a scalar coefficient. Then $p(sxs^{-1}) = p(x)$ for
any invertible $s \in Br_n(m)$, hence $p(\mathbbm{T}') = 0$.
Finally, $p$ acts by 1 on the center of $Br_n(m)$
hence $Z(\BB_n(m)) = p(Z(\BB_n(m)))\subset p(\BB_n(m)) = p(\k T + 
\mathbbm{T}' + [\BB_n(m),\BB_n(m)])
= p(\k T) = \k T$. We thus proved the following.

\begin{prop} For $m \not\in \mathcal{S}_n$ and $m \neq 2$, the Lie
algebra $\BB_n(m)$  is reductive, and its center is spanned by $T$.
\end{prop}

\begin{lemma} \label{lemcentre} Let $\la \in \Irr'_n$. For $m$ outside a finite set of
rational values, $\rho_{\la}(T) \neq 0$.
\end{lemma}

\begin{proof}
Since $\rho_{\la}(T)$ is scalar, $(\dim \la) \rho_{\la}(T) = \tr \rho_{\la}
(T) = \frac{n(n-1)}{2} \tr \rho(t_{12})$. But $\Sp \rho(t_{12}) \subset
\{ -1, 1 , 1-m \}$ and $1-m \in \Sp \rho(t_{12})$ since $\la \in \Irr'_n$.
The conclusion is immediate.
\end{proof}

\subsection{Representations of $Br_n(m)$}

We will need a few technical results on the
representations of $Br_n(m)$ that we gather here. We assume
$m \not\in S_n$.

\begin{figure}
\begin{center}
\resizebox{16cm}{!}{
\includegraphics{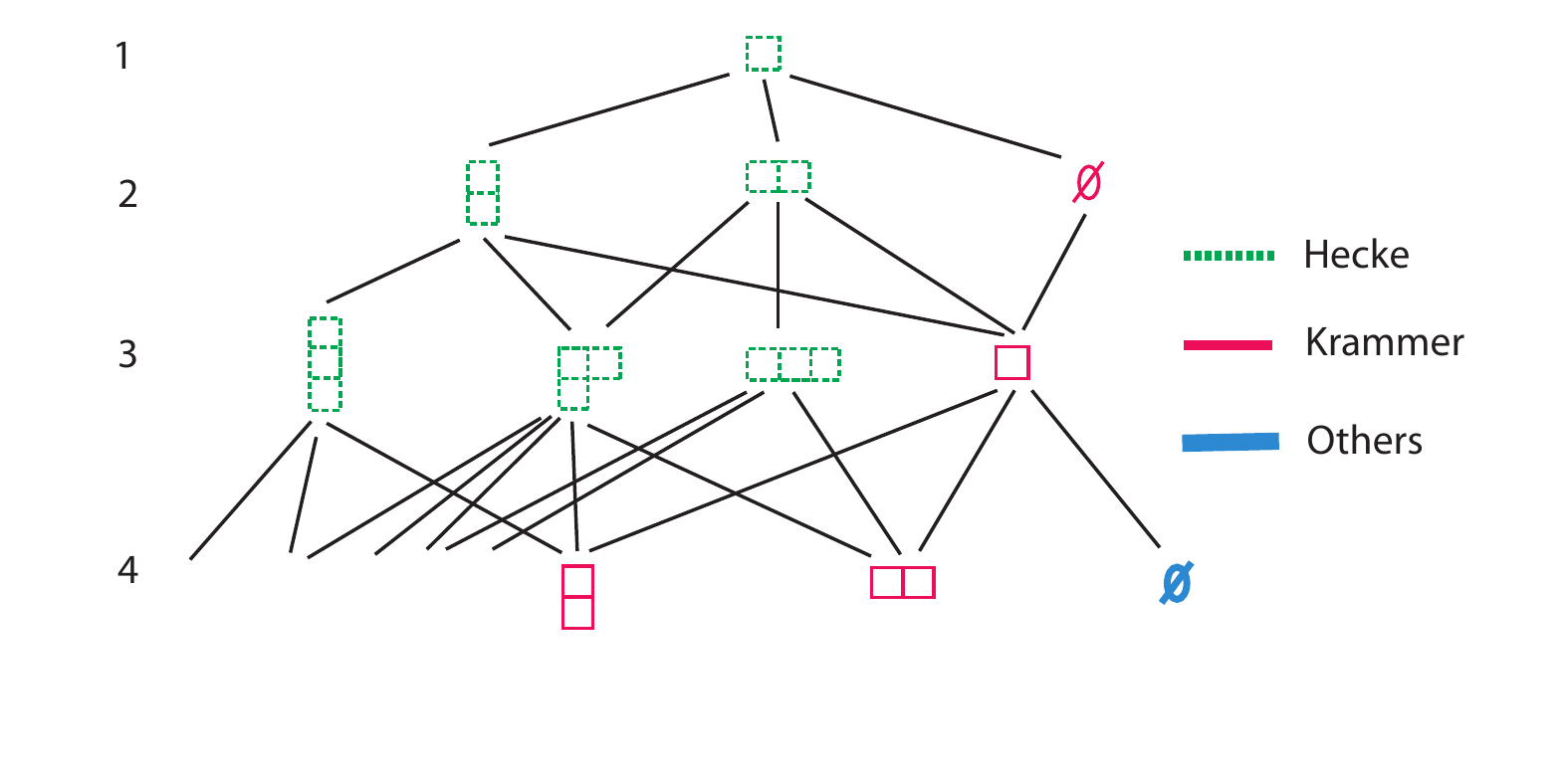}}
\end{center}
\caption{Bratteli diagram for the Brauer algebra}
\label{figbratteli}
\end{figure}

\begin{lemma} \label{lemsp}
Let $n \geq 3$ and $\rho \in \Irr_n$. Either $\Sp \rho(t_{12}) \subset
\{ -1 , 1 \}$ or $\rho \in \Irr'_n$ and $\Sp \rho(t_{12})
= \{ -1, 1 , 1-m \}$.
\end{lemma}
\begin{proof}
If $\rho$ factorizes through $\mathfrak{S}_n$ then clearly
$\Sp \rho(t_{12}) = \Sp \rho(s_{12}) \subset \{ -1,1 \}$.
Otherwise $\rho(p_{12}) \neq 0$, which implies that the
restriction of $\rho$ to $Br_3(m)$ contains $[1]_3$, over
which $t_{12}$ acts with spectrum $\{ -1,1,1-m \}$. This
proves the statement.
\end{proof}

\begin{lemma} Let $n \geq 5$. If $\la \in \Irr(\mathfrak{S}_n)$
and $\dim \la > 1$, then $\dim \la \geq n-1$. If $\la
\in \Irr'_n$ then $\dim \la \geq n(n-1)/2$.
\end{lemma}

\begin{proof} 
The first part is well-known (see e.g. \cite{LIETRANSP} lemme 8) and easy to
check, so we leave it to the reader. We check
the second inequality directly for $n = 5$ and proceed
by induction, assuming $n \geq 6$. We thus assume
$\la \in \Irr'_n$, and denote $\la = (\la_1 \geq \la_2 \geq \dots \geq
\la_r > 0)$ with $|\la| < n$.

We first deal with the cases $|\la| \leq 2$, proving
$\dim [1]_n \geq n(n+1)/2$ and $\dim \la \geq n(n-1)/2$
for $\la \in \{ \emptyset, [2],[1,1] \}$, by induction on $n \geq 5$.
We have $\dim [1]_5 = 15 = 5 \times 6/2$. Let then $n \geq 6$.
If $n$ is even, then $\la \in \{ \emptyset,[2],[1,1] \}$ and the
restriction to $Br_{n-1}(m)$ admits for component $[1]_{n-1}$ ;
it follows that $\dim \la \geq \dim [1]_{n-1} \geq n(n-1)/2$.
If $n$ is odd, then the restriction of $\la = [1]_n$
is $\emptyset_{n-1} + [2]_{n-1} + [1,1]_{n-1}$, hence
$\dim [1]_n \geq 3(n-1)(n-2)/2 \geq n(n+1)/2$ for
$n \geq 10$, and we check that $\dim [1]_9 = 945$ and
$\dim [1]_7 = 105$ also satisfy our assumption.

We now assume $|\la| = 2$. Since $\la \neq \emptyset$ there
exists $\mu_1 \in \Irr'_{n-1}$ such that $\mu_1 \nearrow \la$.
Let $\mu_2 = (\la_1 + 1\geq \la_2 \geq \dots \geq \la_r > 0)_{n-1}$
and $\mu_3 = (\la_1 \geq \la_2 \geq \dots \geq \la_r \geq 1 > 0)_{n-1}$.
We have $\mu_2,\mu_3 \in \Irr_{n-1}$ and $\la \nearrow \mu_2$,
$\la \nearrow \mu_3$. If $\dim \mu_2 = 1$ then
$r = 1$, and in this case $\dim \mu_3 = 1$ implies
$\la_1 = \la_r = 1$, that is $\la = [1]$ which has
been excluded. It follows that $\dim \mu_2 > 1$ or $\dim \mu_3 > 1$,
hence $\dim \mu_2 + \dim \mu_3 \geq 1 + (n-2) = n-1$
(note that $(n-1)(n-2)/2 \geq n-2$), therefore
$\dim \la \geq n-1 + (n-1)(n-2)/2 = n(n-1)/2$, and this concludes
the proof.
\end{proof}

\subsection{Isomorphic representations}

Let $\XX = \{ 1 , -1, m-1 \}$. For every representation
$\rho$, we have $\Sp \rho(t) \subset \XX$. We assume
$m \not\in \mathcal{S}_n$ and $m \neq 2$ in what follows.

Since $\BB_n(m)$ generates $Br_n(m)$ as an algebra, for
every two irreducible representations $\rho^1,\rho^2 \in \Irr_n$,
we have $\rho^1 \simeq \rho^2$ if and only if
$\rho^1_{\BB} \simeq \rho^2_{\BB}$.

For such irreducible representations, we denote $\rho^1_{\BB'}$,
$\rho^2_{\BB'}$ the irreducible representations
of $\BB'_n(m) = [\BB_n(m),\BB_n(m)]$ that they induce. If
$\rho^1 \simeq \rho^2$ then $\rho^1_{\BB'} \simeq \rho^2_{\BB'}$.
We are interested in the converse, and assume that
we have $\rho^1_{\BB'} \simeq \rho^2_{\BB'}$.
We let 
$$
t'_{ij} = t_{ij} - 2T/n(n-1) = \frac{2}{n(n-1)} \sum_{k,l} t_{ij} - t_{kl}.
$$
We saw in \S 5.1 that $p(t'_{ij}) = 0$. Since $p$ is a projector onto
$Z(\BB_n(m))$ whose kernel contains $\BB'_n(m)$, it induces
the
canonical projection of $\BB_n(m)$ onto its center, whence
$t'_{ij} \in \BB'_n(m)$. From this it is clear that $\BB'_n(m)
= [\BB_n(m),\BB_n(m)]$ is generated by the $t'_{ij}$.

Let $t' = t'_{12}$, and assume without loss of generality
that $\rho^1, \rho^2$ share the same underlying vector space $V$.
We have $P \in \GL(V)$ such that $\rho^1(t') = P \rho^2(t') P^{-1}$,
that is $\rho^1(t) = P \rho^2(t) P^{-1} + \alpha$,
with $\alpha = (1/N)(\rho^1(T) - \rho^2(T)) \in \k$ and $N = n(n-1)/2$.
Since $\Sp \rho^i(t) \subset \XX$ we have
$\alpha \in \XX - \XX$, that is $\alpha = 0$ or $\alpha \in 
 \{ 2, 2-m,-m,-2,m-2,m \}$. This latter set
has 6 elements, except when $m \in S = \{ -2, 0, 1, 2 , 4 \}$.
We assume $m \not\in S$. This imposes that, either $\alpha = 0$,
or both $\rho^1(t)$ and $\rho^2(t)$ have a single eigenvalue.
Since the $\rho^i(t)$ are semisimple endomorphisms,
at least for $m \not\in \{ 0, 2 \}$, the $\rho^i(t)$
have then to be scalars. Through conjugation under $\mathfrak{S}_n$
this has for consequence that all the $\rho^i(t_{kl})$ are scalars, 
that is $\dim \rho^1 = \dim \rho^2 = 1$, by
irreducibility of $\rho^1$ and $\rho^2$. This proves the following.

\begin{prop} Let $S = \{ -2, 0, 1, 2 , 4 \}$ and assume $m \not\in \mathcal{S}_n
\cup S$. For $\rho^1,\rho^2 \in \Irr_n$ with $\dim \rho^i > 1$, $\rho^1 \simeq \rho^2$
if and only if $\rho^1_{\BB'} \simeq \rho^2_{\BB'}$.
\end{prop}

We now assume that $\rho^2_{\BB'}$ is isomorphic to the
dual $(\rho^1_{\BB'})^*$ of $\rho^1_{\BB'}$. Let $P \in \GL(V)$
be an intertwinner. We have $- ^t \rho^1(t') = P \rho^2(t')
P^{-1}$ hence $- ^t \rho^1(t) = P \rho^2(t) P^{-1} + \alpha$
with $\alpha = (- \rho^1(T) - \rho^2(T))/N$. Up to conjugation
of the matrix model of $\rho^1$ we can assume that $\rho^1(t)$
is diagonal, with eigenvalues $\la_1,\dots,\la_r \in \XX$.
Then $\alpha - \la_i = \mu_i \in \Sp \rho^2(t) \subset \XX$ and
$\rho^2(t)$ is diagonalizable with eigenvalues $\mu_1,\dots,\mu_r$.
In particular $\alpha \in \XX + \XX$. We write down the
addition table for $\XX$ :
$$
\begin{array}{c||c|c|c|}
& 1 & -1 & m-1 \\
\hline
\hline
1 & 2 & 0 & m \\
\hline
-1 & 0 & -2 & m-2 \\
\hline
m-1 & m & m-2 & 2(m-2) \\
\end{array}
$$
The values in the upper triangular corner of this
table are distinct provided that $m \not\in S^* = \{ 0,2,-2,4,3,1 \}$.
In this case, $\alpha \in \{ 2, -2, 2 (m-2) \}$ implies that
$\forall i \ \la_i = \mu_i = c$ is independent of $i$,
namely that the $\rho^i(t_{kl})$ are scalars and $\dim \rho^i = 1$.
Excluding that case (hence assuming $n \geq 3$), we thus have $\alpha \in \{0,m,m-2 \}$,
which corresponds to $\Sp \rho^i (t) \in \{ \{ 1,-1 \},
\{1, m-1 \}, \{ -1 , m-1 \} \}$. By lemma \ref{lemsp}, we have
$\Sp \rho (t) \not\in \{ 
\{1, m-1 \}, \{ -1 , m-1 \} \}$ (when $m -1 \not\in \{ 1, -1 \})$.
It follows that $\Sp \rho^1(t) = \Sp \rho^2(t) = \{ 1, -1 \}$,
meaning that the $\rho^i$ factorize through $\k \mathfrak{S}_n$,
which is the quotient of $Br_n(m)$ by $p_{12}$
(from $\rho(t)^2 = 1$ and $\rho(t^2 - 1) = (m-2) \rho(p_{12})$
we indeed get $\rho(p_{12}) = 0$ for $m \neq 2$). But
this case is known by \cite{LIETRANSP} (proposition 2 and \S 5.1 ; see also
\cite{IH2} proposition 2.7), so this concludes the
proof of the following proposition :

\begin{prop} Let $S^* = \{ -2, 0, 1, 2 ,3, 4 \}$ and assume
$m \not\in \mathcal{S}_n
\cup S^*$. For $\rho^1,\rho^2 \in \Irr_n$, $\rho^2_{\BB'} \simeq
(\rho^1_{\BB'})^*$
if and only if, either $\dim \rho^1 = \dim \rho^2 = 1$, or
$\rho^1$ and $\rho^2$ factorize through $\k \mathfrak{S}_n$
and $\rho^2 = \rho^1 \otimes \eps$, where $\eps$ is the sign character
of $\mathfrak{S}_n$.
\end{prop}

\subsection{Infinitesimal Hecke algebra}

Let $\mathcal{H}_n$ denote the Lie subalgebra of $\k \mathfrak{S}_n$
generated by the transpositions.
This Lie subalgebra, which is a special case of
the infinitesimal Hecke algebras dealt with in \cite{IH2}, is reductive
and has been decomposed in \cite{LIETRANSP}. The quotient map $Br_n(m)
\to \k \mathfrak{S}_n$
induces a Lie algebra morphism $\BB_n(m) \to \mathcal{H}_n$, that
is onto and maps center to center isomorphically. It thus
induces $\BB'_n(m) \onto \mathcal{H}_n'$. Let then
$\Irr'_n = \{ \la \in \Irr_n \ | \ |\la| < n \}$,
where $|\la|$ denotes the number of boxes of the Young diagram
associated to $\la$. We still assume $m \not\in \mathcal{S}_n$
and consider the isomorphism
$Br_n(m)'  \to \bigoplus_{\la \in \Irr_n} \sl(V_{\la})$, where $V_{\la}$
is the underlying vector space of $\la$. Since $\la \in \Irr_n \setminus
\Irr'_n$ means that $\la$ factorizes through $\k \mathfrak{S}_n$,
this isomorphism can be written $Br_n(m)' \to (\k \mathfrak{S}_n)' \times
\bigoplus_{\la \in \Irr'_n} \sl(V_{\la})$, hence induces an injective
map
$$
\BB'_n(m) \into \mathcal{H}_n' \times \bigoplus_{\la \in \Irr'_n} \sl(V_{\la})
$$
Our goal here is to show that this map is surjective, namely
\begin{theo} \label{theobrauer} Let $m \not\in \mathcal{S}_n \cup S \cup S^* $. Then $\BB_n(m)$ is
a reductive algebra with 1-dimensional center, whose
derived Lie algebra is naturally isomorphic to
$\mathcal{H}_n' \times \bigoplus_{\la \in \Irr'_n} \sl(V_{\la})$.
\end{theo}

For the convenience of the reader,
we recall from \cite{LIETRANSP} the simple ideals
of $\mathcal{H}'_n$. For $\rho \in \Irr(\mathfrak{S}_n)$, identified
with a partition of $n$ or a Young diagram of size
$n$, we denote $V_{\rho}$ the underlying vector space and
$\rho_{\mathcal{H}'} : \mathcal{H}'_n \to \sl(V_{\rho})$ the
induced representation of $\mathcal{H}'_n$. Then the orthogonal
$\mathcal{H}(\rho)$ of $\Ker \rho_{\mathcal{H}'}$ for
the Killing form is a simple ideal of $\mathcal{H}'_n$, and
all simple ideals are obtained this way.
A non-overlapping list is given by $\mathcal{H}([n-1,1]) \simeq
\sl_{n-1}(\C)$ and, for $\rho$ \emph{not a hook}, letting
$\eps : \mathfrak{S}_n \onto \{ \pm 1 \}$ denote the sign character,
\begin{itemize}
\item $\mathcal{H}(\rho) \simeq \so(V_{\rho})$ when the symmetric
square $S^2\rho$ contains $\eps$ ;
\item $\mathcal{H}(\rho) \simeq \sp(V_{\rho})$ when the symmetric
square $S^2\rho$ contains $\eps$ ;
\item $\mathcal{H}(\rho) = \mathcal{H}(\rho \otimes \eps) \simeq \sl(V_{\rho})$
otherwise.
\end{itemize}

\bigskip

For a given $n$, the decomposition of $\BB'_n(m)$ given by the theorem
is equivalent to the property
that $\rho_{\la}(\BB'_n(m)) = \sl(V_{\la})$ for all
$\la \in \Irr'_n$. Indeed, if we have this property, then
the semisimple Lie algebra $\BB'_n(m)$ contains simple Lie
ideals isomorphic to $\sl(V_{\la})$ for every $\la \in \Irr_n'$,
which do not intersect the simple Lie ideals inherited from
$\mathcal{H}_n'$, and do not intersect each other : indeed,
if there was such an intersection this would mean that
this two simple ideals coincide, meaning that
two $\la,\mu \in \Irr_n'$ of the same dimension would factorize
through the same ideal. Since $\sl_N(\C)$ admits at most two
irreducible representations of dimension $N$, this would
imply $\la_{\BB'} \simeq \mu_{\BB'}$ or
$\la_{\BB'} \simeq (\mu_{\BB'})^*$, which is excluded by the propositions
above ; the case $\la \in \Irr'_n, \mu \in \Irr_n \setminus \Irr'_n$ is dealt
with similarly. This property for any given $n$ thus implies the theorem for
this $n$.

\subsection{Induction step}

We prove here that $\rho(\BB'_n(m)) = \sl(V_{\la})$ for
every $\la \in \Irr'_n$ by induction on $n$, assuming
the theorem true for $n-1$. Note that no confusion should arise in the notations
between $\Irr_n$ and $\Irr_{n-1}$ as, if $\la \in \Irr_n$ and $\mu \in \Irr'_n$,
then $|\la|$ and $|\mu$ does not have the same parity.
We let $V_{\la}$ denote the underlying space of $\la \in \Irr_n$ or
$\la \in \Irr_{n-1}$. For two partitions $\la,\mu$, we use the notation
$\la \nearrow \mu$ for $|\mu| = |\la|+1$ and $\mu_i = \la_i$ for all but
one $i$. When confusion could arise, we let $\rho_{\la}$ denote
the representation of $Br_n(m)$ or $Br_{n-1}(m)$ associated to the
partition $\la$.

The branching
rule for $\BB_{n-1}(m) \subset \BB_n(m)$ can be written
as
$$
\Res \rho_{\la} = \sum_{\mu \nearrow \la } \rho_{\mu} + \sum_{\la
\nearrow \mu} \rho_{\mu}
$$
We consider several cases. First note that
we can assume $n \geq 5$. Indeed, for $n \leq 4$,
an element of $\Irr'_n$ is either
the infinitesimal Krammer
representation
(or its transformed under the automorphism $s_{ij} \mapsto -s_{ij}$,
$p_{ij} \mapsto -p_{ij}$ of $Br_n(m)$)
, which has been studied separately in \cite{KRAMMINF}, or $[\emptyset]_4$.
In this last case, its restriction to $Br_3(m)$ is the irreducible
representation $[1]_3$ which is a Krammer representation, hence
$\rho_{[\emptyset]_4}(\BB'_4(m)) \supset \rho_{[1]_3}(\BB'_3(m)) 
= \sl(V_{[1]_3})$.

\subsubsection{$|\la| < n-2$}
In this case, for every $\mu \nearrow \la$ or $\la \nearrow \mu$,
we have $\mu \in \Irr'_{n-1}$. By the induction assumption,
$\rho_{\la}(\BB'_n(m))$ contains
$$
\left( \bigoplus_{\mu \nearrow \la} \sl(V_{\mu}) \right)
\oplus
\left( \bigoplus_{\la \nearrow \mu} \sl(V_{\mu}) \right)
$$
which has (semisimple) rank
$$
\sum_{\mu \nearrow \la} (\dim V_{\mu} - 1) + \sum_{\la
 \nearrow  \mu} (\dim V_{\mu} - 1).
$$
For $\mu \in \Irr'_{n-1}$ and $n-1 \geq 3$, that is $n \geq 4$, we
have $\dim V_{\mu} \geq 3$ hence $\dim V_{\mu} > (1/2) \dim V_{\mu}$
hence $\rk \rho_{\la}(\BB'_n(m)) > (1/2) \dim V_{\la}$.
This implies $\rho_{\la}(\BB'_n(m)) = \sl(V_{\la})$
(see \cite{KRAMMINF} prop. 3.8).

\subsubsection{$|\la| = n-2$ and no $\mu$ with $\la \nearrow \mu$ is a hook}

In this section, the original assumption $n \geq 5$ implies $n \geq 6$,
because for $n = 5$ and $|\la| = 3$ there is always a hook
$\mu$ with $\la \nearrow \mu$.

First notice that $\la \nearrow \mu$ and $\la \nearrow \mu'$,
where $\mu'$ is the transposed partition of $\mu$, implies
that either $\mu = \mu'$, in which case $\mu$ is uniquely determined,
or $\la = \la'$ and $\la \nearrow \nu \Rightarrow \la  \nearrow \nu'$
for every $\nu$.
If there is no such $\mu$, for $n \geq 6$ we get again
$\dim V_{\mu} \geq 3$ for all irreducible component $\mu$ of the
restriction (as $\dim V_{\mu} = 1$ would
imply that $\la$ is a hook, and if $\la$ is a hook there
is a hook $\mu$ with $\la \nearrow \mu$, contradicting the assumption),
and we conclude as in the previous case that
$\rho_{\la}(\BB'_n(m)) = \sl(V_{\la})$. 

If there is such a $\mu$ with $\mu = \mu'$, it
is uniquely determined, we have
$\rk \rho_{\mu}(\BB_{n-1}(m)') = (\dim V_{\mu})/2$, and the same
computation as before shows again $\rk \rho_{\la}(\BB'_n(m)) > (1/2) \dim V_{\la}$.

So the only case that we have to deal with is $\la = \la'$. Then $\la
\nearrow \mu \Rightarrow \mu \neq \mu'$ and $\la \nearrow \mu'$.
Moreover,
$$
\rk \rho_{\la}(\BB_{n-1}(m)') = \sum_{\mu \nearrow \la} (\dim V_{\mu} -1)
 + \frac{1}{2} \sum_{\la \nearrow \mu} (\dim V_{\mu} - 1).
$$
If $\dim V_{\mu} \geq 3$ (always the case for $n \geq 6$), we have
$(1/2)(\dim V_{\mu} - 1) \geq (1/3) \dim V_{\mu}$ hence
$$
r = \rk \rho_{\la}(\BB_{n-1}(m)') > \frac{1}{3} \sum_{\mu \nearrow \la} \dim
V_{\mu} + \sum_{\la \nearrow \mu} \dim V_{\mu} ) = \frac{1}{3} \dim V_{\la} = d/3.
$$
Then $r > d/3$ implies $d < 3r < (r+1)^2$, because $d \geq 3$.
This implies that $\rho_{\la}(\BB')$
is simple (see \cite{IH2} lemma 3.3 (I)). Moreover $\rk \rho_{\la}(\BB') > (\dim V_{\la})/3 >
(\dim V_{\la})/4$. Note that, for every $\mu \in \Irr'_{n-1}$,
in particular for every $\mu \nearrow \la$,
we have $\dim V_{\mu} \geq (n-1)(n-2)/2$. If $|\la| \geq 1$
we thus get $\rk \rho_{\la}(\BB') \geq (n-1)(n-2)/2 -1 \geq 9$
for $n \geq 6$. These two inequalities imply $\rho_{\la}(\BB') =
\sl(V_{\la})$ (see \cite{IH2} lemma 3.4). If $|\la| = 0$,
then $\rho_{\la}(\BB') = \rho_{[1]_{n-1}}(\BB') = \sl(V_{[1]_{n-1}})
= \sl(V_{\la})$ by the induction assumption.

\subsubsection{$|\la| = n-2$ and there exists a hook $\mu$ with $\la \nearrow
\mu$}

In that case $\la$ has the shape of a hook $[n-k,1^{k-2}]$ and we
can assume that $k \geq 3$, $n-k-2 \geq 2$ (and $n \geq 5$). Indeed,
it is otherwise isomorphic to the infinitesimal Krammer
representation
(or its transformed under the automorphism $s_{ij} \mapsto -s_{ij}$,
$p_{ij} \mapsto -p_{ij}$ of $Br_n(m)$).

We denote $C(n,k)$ the dimension of this representation. The $\mu$ such that
$\la \nearrow \mu$ are the partitions $A = [n-k,1^{k-1}]$, $B = [n-k+1,1^{k-2}]$,
$C = [n-k,2,1^{k-3}]$. The $\mu$ such that $\mu \nearrow \la$
are the partitions $E = [n-k-1,1^{k-2}]$, $F = [n-k,1^{k-3}]$.
We have
$$
\dim A = \cnp{n-2}{k-1}, \dim B = \cnp{n-2}{k-2}
$$
and, by \cite{LIETRANSP} \S 6.3,
$$
 \dim C = \frac{k-2}{n-k}
\cnp{n-3}{n-k-2} (n-1) = \cnp{n-1}{k-1} \frac{(n-k-1)(k-2)}{n-2}.
$$
Moreover, $\dim E = C(n-1,k)$ and $\dim F = C(n-1,k-1)$.
We have
$$
C(n,k) = \dim V_{\la} = C(n-1,k) + C(n-1,k-1) + \cnp{n-2}{k-1} + 
\cnp{n-2}{k-2} + \dim C
$$
hence $C(n,k) > \cnp{n-1}{k-1} + C(n-1,k) + C(n-1,k-1)$.

We show by induction on $n$ that $C(n,k) > (n-2) \cnp{n-1}{k-1}$.
This holds true for $n = 5$, because $C(5,3) = \dim [2,1]_5 = 20$.
Then
$$
\cnp{n-1}{k-1} + C(n-1,k) + C(n-1,k-1) > \cnp{n-1}{k-1}
+ (n-3) \left( \cnp{n-2}{k-1} + \cnp{n-2}{k-2} \right)
$$
and the last term is equal to $(n-2) \cnp{n-1}{k-1}$,
which proves by induction the inequality, except for $k = 3$. In that
case,
$C(n,3) = \cnp{n-1}{2} + C(n-1,3) + \dim [n-3]_{n-1} + \dim
[n-3,2]_{n-1}$ and $\dim [n-3]_{n-1} = \cnp{n-1}{2}$. Then, by
the induction assumption,
$$
C(n,3) > \cnp{n-1}{2} + (n-3) \cnp{n-2}{2} + \cnp{n-1}{2}
+ \frac{(n-1)(n-4)}{2}
$$
hence
$$
C(n,3) > \frac{(n-1)(n-2)}{2} \left( n-2 + 1 - \frac{n^2-5n+8}{n^2-3n+2}
\right) > \frac{(n-1)(n-2)}{2} (n-2).
$$
In particular, we have
$$
\begin{array}{lclclcl}
\dim E &=& C(n-1,k) &>& (n-3) \cnp{n-2}{k-1} &=& (n-3) \dim A \\
\dim F &=& C(n-1,k-1) &>& (n-3) \cnp{n-2}{k-2} &=& (n-3) \dim B 
\end{array}
$$
if $k \geq 4$ ($\dim B = n-2$ if $k = 3$),
and $\dim E + \dim F > (n-3) \cnp{n-1}{k-1} \geq 4 \dim C$
if and only if $(n-3)(n-2) \geq 4 (n-k-1)(k-2)$, which holds true.
As a consequence, for $k \geq 4$,
$$
\begin{array}{lcl}
\rk \rho_{\la}(\BB') &\geq &\dim E -1 + \dim F - 1 + (n-2) + (\dim C)/2 \\
& \geq & \dim E + \dim F - 2 + (n-2) + (\dim C)/2.
\end{array}
$$
Moreover, $\dim V_{\la} = \dim E + \dim F + \dim A + \dim B
+ \dim C < (\dim E + \dim F) (1 + \frac{1}{n-3} + \frac{1}{4} )$
hence
$
\dim V_{\la} < \frac{5n-11}{4n-12} ( \dim E + \dim F)$.
We thus have
$$
\rk \rho_{\la}(\BB') \geqslant \dim E + \dim F > \frac{4n-12}{5n-11}
\dim V_{\la} > \frac{1}{2} \dim V_{\la}
$$
for $n \geq 5$, hence $\rho_{\la}(\BB') = \sl(V_{\la})$
for $n \geq 5$ and $k \geq 4$.

The only remaining case is when $k = 3$, i.e. $\la$
corresponds to the partition $[n-3,1]$. We still have
$\dim E > (n-3) \dim A$, but this time $\dim F = (n-1)(n-2)/2$,
$\dim B = n-2$,
$$\dim C = \cnp{n-1}{2} \frac{n-4}{n-2} = \frac{(n-1)(n-4)}{2}.
$$
We have $\dim A = (n-2)(n-3)/2$ and
$\dim E > (n-3) \dim A = (n-2)(n-3)^2/2$. Since $n-3 \geq (n-1)/2$
when $n \geq 5$, we also have $\dim E > \frac{n-1}{2} \dim A$. From
$\dim F = \frac{n-1}{2} \dim B$ we then get
$$
\dim E + \dim F > \frac{n-1}{2} \left( \dim A + \dim B \right)
$$
hence
$$
\dim C < \frac{n-1}{n-2} \frac{n-4}{n-3} \frac{1}{n-3} \left(
\dim E + \dim F \right) < \frac{1}{n-3} \left(
\dim E + \dim F \right).
$$
We thus have
$$
\dim V_{\la} < (\dim E + \dim F) \left( 1 + \frac{2}{n-1} + \frac{1}{n-3}
\right) = (\dim E + \dim F) \frac{n^2-n-4}{(n-1)(n-3)}
$$
and $\rk \rho_{\la}(\BB') \geq \dim E + \dim F + \dim C - 3 + (n-2)$
if $n \geq 6$, $\rk \rho_{\la}(\BB') \geq \dim E + \dim F + \frac{\dim C}{2} - 2 + (n-2)$
if $n = 5$ (in that case $C' = C$). In both cases
$$
\rk \rho_{\la}(\BB') > \dim E + \dim F > \frac{(n-1)(n-3)}{n^2-n-4}
\dim V_{\la}.
$$
Finally, $(n-1)(n-3)/(n^2-n-4) \geq 1/2$ for $n \geq 5$, and
this again proves $\rho_{\la}(\BB') = \sl(V_{\la})$ by
the same criterium (\cite{LIETRANSP} proposition 3.8).

\section{Proof of theorem \ref{maintheo}}

We first let $K = \C((h))$. For $\la \in \Irr_n$, $m \in \C$,
we denote $V_{\la}$ the underlying vector space of $\la : Br_n(m)
\to \End(V_{\la})$ over $\C$, and $R_{\la}^{(m)} :
B_n \to \GL(V_{\la} \otimes K)$ the representation constructed in
\S 4 from $\rho$ by using some associator. Let $\Gamma_{\la}$
denote the algebraic (Zariski) closure of $R_{\la}^{(m)}(B_n)$
inside $\GL(V_{\la} \otimes K)$.

Assume $\la \in \Irr'_n$ and $m \not\in \mathcal{S}_n \cup S \cup S^*$.
By theorem \ref{theobrauer} and proposition \ref{irred}, the Lie algebra
$\mathrm{Lie}\, \Gamma_{\la}$
of $\Gamma_{\la}$ contains $\sl(V_{\la}) \otimes K$. Recall
from the proof of lemma \ref{lemcentre} that $\la(T)$ is a (scalar) rational affine polynomial
in $m$ ; we denote $E_{\la,n}$ the set of rationals $m$ such that
$\la(T) = 0$. The generator $\gamma_n = (\sigma_1 \dots \sigma_{n-1})^n$
of $Z(B_n)$ is mapped to $\exp T$ through
$B_n \to \mathfrak{S}_n \ltimes \exp \widehat{\mathcal{T}_n}$ ;
thus $\mathrm{Lie}\, \Gamma_{\la}$ contains $K$ as soon as $m \not\in
E_{\la,n}$.

It follows that $\Gamma_{\la} = \GL(V_{\la} \otimes K)$ if $m \not\in E_{
\la,n}$, and $\Gamma_{\la} \supset \SL(V_{\la} \otimes K)$. In
particular the algebraic closure of $R_{\la}^{(m)}(B_n')$
is equal to $\SL(V_{\la} \otimes K)$.

We now consider $BMW_n(e^h,e^{(1-m)h})$ as an algebra over $K$, and still
assume $m \not\in \mathcal{S}_n \cup S \cup S^*$. Then $BMW_n(e^h,
e^{(1-m)h})$ is the direct sum of the Hecke algebra $H_n(e^h)$ and
of $\bigoplus_{\la \in \Irr'_n} \End(V_{\la} \otimes K)$. Let
$G_0$ denote the algebraic closure of $B_n'$ inside $H_n(e^h)$. This
closure has been descrived in full detail in \cite{LIETRANSP}. We recall
from there
that its
Lie algebra is $\mathcal{H}_n' \otimes K$, and that it is a direct
sum of $\SL_N,\SP_N$ and $\SO_N$ with respect to some
hyperbolic nondegenerate quadratic form.
From theorem \ref{theobrauer}
we get similarly that the image of $B_n'$ inside $BMW_n(e^h,e^{(1-m)h})$
is Zariski-dense inside $G_0 \times \prod_{\la \in \Irr'_n}
\SL(V_{\la} \otimes K)$, whose Lie algebra is $\BB'_n(m)$.

When $m \not\in \Q$, $e^h$ and $e^{(1-m)h}$
are algebraically independent over $\Q$, so we get an embedding $\Q(s,\alpha)
\into K$ through $s \mapsto e^h$ and $\alpha \mapsto e^{mh}$.
Similarly, embedding $\Q(s)$ into $K$ through $s \mapsto e^h$
we get realizations of the representations of $BMW_n(s,s^m)
\otimes K$ for $m \in \Z$ and $m \not\in \mathcal{S}_n \cup S \cup S^*$,
where $BMW_n(s,s^m)$ is defined over $\Q(s)$ through the $R_{\la}^{(m)}$,
$\la \in \Irr_n$. We know the Zariski closure of the $R_{\la} = R_{\la}^{(m)}(B_n')$
for $\la \not\in \Irr'_n$ from \cite{LIETRANSP}, provided we know that
the orthosymplectic groups involved there when $\la = \la'$ are defined
over $\Q(s)$. The bilinear form defining them span
the subspace of $R_{\la} \otimes R_{\la}$ over which $B_n$ acts
by the sign character $B_n \onto \{ \pm 1 \}$, $\sigma_i \mapsto -1$.
Since $\eps$ and the $R_{\la}$ are defined over $\Q(s)$, this
subspace has nonzero points over $\Q(s)$, so these groups
are indeed defined over $\Q(s)$.


Theorem \ref{maintheo} is then an immediate
consequence of the above for an arbitrary (characteristic 0) field,
by noticing that all the algebraic groups $G$ involved here are defined
over $\Q(s)$, satisfy that $G(L_1)$ is Zariski-dense in
$G(L_2)$ for $\Q(s) \subset G(L_1) \subset G(L_2)$, and browsing
along the following pattern of field extensions.

\begin{center}
\mbox{\xymatrix{
 \C((h)) & & K \\
 & \ar@{-}[ul] \ar@{-}[ur] \ar@{-}[d] \Q(s,\alpha) & \\
 & \ar@{-}[uul] \ar@{-}[uur]\Q(s) & 
}}
\end{center}

One point that remains to be clarified is the type of the
orthogonal groups possibly appearing for $\la \not\in \Irr'_n$
and $\la = \la'$, when considered over $\Q(s)$. We know that
the quadratic forms $\beta_{\la}$ involved here are hyperbolic over
$\C((h))$, but it was not proved in \cite{LIETRANSP} that
they are hyperbolic over $\Q(s)$. We prove this below.

\begin{lemma} 
For $n \geq 2$, $\la \vdash n$, $\la = \la'$,
if $\eps \into S^2 R_{\la}$, then $\beta_{\la}$ is hyperbolic.
\end{lemma}
\begin{proof} We prove the lemma by induction on $n$, the
cases $n = 2$ being clear. By Young rule the restriction
to $B_{n-1}$ of $R_{\la}$ is the direct sum of the $R_{\mu}$
for $\mu \nearrow \la$. Notice that $\mu \nearrow \la$
implies $\mu' \nearrow \la$ under our assumptions. Let $\mu,\nu
\nearrow \la$ and consider the restriction $\beta : V_{\mu} \otimes
V_{\nu} \to \Q(s)$ of $\beta_{\la}$. First assume $\nu \neq \mu'$.
If $\beta \neq 0$ it would provide an isomorphism $R_{\mu}
\simeq R_{\nu}^* \otimes \eps \simeq R_{\nu'}$
where $\eps$ is the sign character of $B_{n-1}$ and $R_{\nu}^*$ is the
dual representation of $R_{\nu}$. Since $\mu \neq \nu'$
we have $\beta = 0$. If $\nu = \mu' = \mu$, by induction we
know that $\beta$ is either 0 or hyperbolic, but the case
$\beta = 0$ is excluded because $\beta_{\la}$ would
then be degenerate. If $\nu = \mu' \neq \mu$, we consider the
restriction of $\beta_{\la}$ to $V_{\mu} \oplus V_{\mu'}$. If
$M : V_{\mu'} \to V_{\mu}$ affords $R_{\mu'} \simeq R_{\mu}^+ \otimes \eps$,
this restriction can be written in matrix form as
$\begin{pmatrix} 0 &  M \\ ^t M & 0 \end{pmatrix}$ and is
therefore hyperbolic. The quadratic space $(V_{\la}, \beta_{\la})$
is thus isomorphic to a direct sum of hyperbolic spaces,
and is therefore hyperbolic.
\end{proof}

\tableofcontents



\end{document}